\documentclass{amsart}
\usepackage[utf8]{inputenc}
\usepackage{amsthm}
\usepackage{amsmath}
\usepackage{graphicx}
\usepackage{hyperref}
\usepackage{mathrsfs}

\usepackage{mathtools}

\title{Rigidity of diagonally embedded triangle groups}
\author{Jean-Philippe Burelle}
\address{CNRS \and Institut des Hautes \'Etudes Scientifiques\\
35 route de Chartres, 91440 Bures-sur-Yvette, France}
\email{jburelle@ihes.fr}
\thanks{
This project has received funding from the European Research Council (ERC) under the European Union’s Horizon 2020 research and innovation programme (ERC starting grant DiGGeS, grant agreement No 715982).}
\date{May 2019}

\newtheorem{prop}{Proposition}[section]
\newtheorem{lem}{Lemma}[prop]
\newtheorem{thm}{Theorem}[prop]
\newtheorem{thmintro}{Theorem}
\newtheorem{cor}{Corollary}[prop]

\theoremstyle{definition}
\newtheorem{defn}{Definition}[prop]
\newtheorem{nota}{Notation}[prop]

\newcommand{\Refl}[2]{\mathsf{R}_{#1}^{#2}}
\newcommand{\GL}{\mathsf{GL}}
\newcommand{\Gr}{\mathsf{Gr}}
\newcommand{\M}{\mathsf{M}}
\newcommand{\Sp}{\mathsf{Sp}}
\newcommand{\PSp}{\mathsf{PSp}}
\newcommand{\Sppm}{\mathsf{Sp}^\pm}
\newcommand{\graph}{\mathsf{Graph}}
\newcommand{\bR}{\mathbb{R}}

\newcommand{\bZ}{\mathbb{Z}}
\newcommand{\bH}{\mathbb{H}}
\newcommand{\RP}{\mathbb{RP}}
\newcommand{\SL}{\mathsf{SL}}
\newcommand{\PSL}{\mathsf{PSL}}
\newcommand{\PGL}{\mathsf{PGL}}
\newcommand{\Lag}{\mathsf{Lag}}

\newcommand{\Hom}{\mathsf{Hom}}
\newcommand{\Inn}{\mathsf{Inn}}
\newcommand{\conf}{\mathsf{Conf}}
\newcommand{\rgeom}{\rho_\mathrm{geom}}
\newcommand{\rirr}{\rho_\mathrm{irr}}

\DeclareMathOperator{\Ort}{O}
\DeclareMathOperator{\Aut}{Aut}
\DeclareMathOperator{\stab}{Stab}
\newcommand{\Span}{\mathsf{span}}

\begin{document}

\begin{abstract}
We show local rigidity of hyperbolic triangle groups generated by reflections in pairs of $n$-dimensional subspaces of $\bR^{2n}$ obtained by composition of the geometric representation in $\PGL(2,\bR)$ with the diagonal embeddings into $\PGL(2n,\bR)$ and $\PSp^\pm(2n,\bR)$.
\end{abstract}

\maketitle
\tableofcontents

\section{Introduction}
% To mention in intro : Rigidity of lattices. Rigidity/flexibility of ideal triangle groups in complex hyperbolic spaces. Maximal representations (diagonal) of surface groups. 

We investigate subgroups of $\PGL(2n,\bR)$ generated by reflections in pairs $U,V$ of half-dimensional subspaces and their deformation properties. Such a reflection is an element of order two in $\PGL(2n,\bR)$ which has eigenvalues $\pm 1$, each with multiplicity $n$.

Denote by
  \[\Delta(k_1,k_2,k_3) = \langle a,b,c ~|~ a^2=b^2=c^2=(ab)^{k_1}=(bc)^{k_2}=(ac)^{k_3}=1 \rangle\]
the triangle group with parameters $k_1,k_2,k_3$. A triangle group is called \emph{hyperbolic} if it can be realized as a discrete subgroup of isometries of the hyperbolic plane generated by reflections in the sides of a triangle. This is equivalent to the inequality $\frac{1}{k_1}+\frac{1}{k_2} + \frac{1}{k_3} < 1$. The geometric representation $\rgeom : \Delta(k_1,k_2,k_3) \rightarrow \PGL(2,\bR)$ obtained by mapping the generators $a,b,c$ to the reflections in the sides of a hyperbolic triangle with interior angles $\frac{\pi}{k_1},\frac{\pi}{k_2},\frac{\pi}{k_3}$ is an example of a group generated by reflections in pairs of half-dimensional subspaces of $\bR^2$.

One way to obtain a triangle group which is generated by reflections in half-dimensional subspaces of $\bR^{2n}$ is by postcomposition with the irreducible representation $\rirr : \PGL(2,\bR) \rightarrow \PGL(2n,\bR)$. In \cite{longth}, Long and Thistlethwaite compute the dimension of the $\PGL(m,\bR)$-\emph{Hitchin component}, the connected component of representations of a triangle group containing $\rirr\circ \rho_{\mathsf{geom}}$. The analog for the groups $\Sp^\pm(2n,\bR)$ was computed by Weir \cite{weir}. Except in a few cases of small dimension, these components have positive dimension.

In their recent paper \cite{ALFHitchin}, Alessandrini, Lee and Schaffhauser initiate the study of \emph{higher Teichm\"uller spaces} for all orbifold surfaces. They prove that the Hitchin component of any orbifold is homeomorphic to a ball, and extend the dimension counts of Long and Thistlethwaite to all orbifold surface groups.

In this paper, we will be interested in representations of triangle groups which factor through the diagonal embedding of $\PGL(2,\bR)$ into $\PGL(2n,\bR)$. By contrast with the Hitchin case, we obtain the following rigidity result :

\begin{thmintro}
  Let  $\Gamma = \Delta(k_1,k_2,k_3)$ be a hyperbolic triangle group and $n\geq1$. Let $i_{2n}$ denote the diagonal representation $\PGL(2,\bR) \rightarrow \PGL(2n,\bR)$.
  Then, any continuous deformation of the composition $i_{2n}\circ\rgeom$ is conjugate to $i_{2n}\circ\rgeom$. 
\end{thmintro}

Alessandrini, Lee and Schaffhauser introduce a notion of \emph{expected dimension} of a component of the character variety for an orbifold surface group. They show that this count corresponds to the actual dimension in the case of the Hitchin components. For the component containing the diagonal representation as in Theorem 1, this expected dimension count is negative.

They also obtain results about Hitchin representations into $\PSp^\pm(2n,\bR)$. They prove that except for $k_1,k_2\leq 3$ in $\PSp^\pm(4,\bR)$, Hitchin representations of triangle groups into $\PSp^\pm(2n,\bR)$ admit non-trivial deformations.

Using the same techniques as for $\PGL(2n,\bR)$, we find many examples which are rigid in arbitrary symplectic groups. This is because there are many non-conjugate diagonal representations of $\PGL(2,\bR)$ into $\Sp^\pm(2n,\bR)$. These diagonal representations are parameterized by the possible signatures $(p,q)$ of a nondegenerate symmetric bilinear form on $\bR^n$ up to exchanging $p$ and $q$.

\begin{thmintro}\label{thm:sympl}
  Let  $\Gamma = \Delta(k_1,k_2,k_3)$ be a hyperbolic triangle group.
  Then, for any $(p,q)$-diagonal representation $i_{(p,q)}$, any deformation of the composition $i_{(p,q)}\circ\rgeom$ is conjugate to $i_{(p,q)}\circ\rgeom$.
  In particular, there are at least $\lceil\frac{n+1}{2}\rceil$ isolated points in the character variety
  \[\chi\left(\Gamma,\Sp^\pm(2n,\bR)\right) = \Hom(\Gamma,\Sp^\pm(2n,\bR))/\Inn(\Sp^\pm(2n,\bR)).\]
\end{thmintro}

We obtain a stronger result for the positive-definite case :
\begin{thmintro}\label{thm:symplpositive}
  Let  $\Gamma = \Delta(k_1,k_2,k_3)$ be a hyperbolic triangle group. Let $\rho:\Gamma\rightarrow \PSp^\pm(2n,\bR)$ be a representation mapping the generators $a,b,c$ of $\Gamma$ to reflections in pairs of Lagrangian subspaces. If there are three eigenspaces $L_a,L_b,L_c$ of the images $\rho(a),\rho(b),\rho(c)$ such that the Maslov form of the triple is positive-definite, then $\rho$ is locally rigid.
\end{thmintro}

For the fundamental group of a closed surface, the composition of a discrete and faithful representation into $\PGL(2,\bR)$ with the positive definite diagonal embedding $i_{(n,0)}$ into $\PSp^\pm(2n,\bR)$ gives rise to a maximal representation, that is, a representation with maximal Toledo invariant. Spaces of maximal representations and Hitchin components are examples of \emph{higher Teichm\"uller theories}, connected components of the space of representations which consist entirely of discrete and faithful representations. For a recent survey on higher Teichm\"uller theories, see \cite{HTTSurvey}.

As a corollary of Theorem \ref{thm:sympl}, for a hyperbolic triangle group there are no non-Fuchsian representations in the \emph{diagonal component} of maximal representations into $\PSp^\pm(2n,\bR)$. This is in contrast to the surface group case where the analogous component always contains Zariski-dense representations \cite{BIW,gwcomponents}.

The strategy of proof will be the use an endomorphism-valued invariant of quadruples of subspaces generalizing the cross-ratio of four points in $\RP^1$. For diagonally embedded representations, these invariants are scalar multiples of the identity. This fact, together with the rigidity of finite order elements, is what allows us to conclude that the representations are rigid.

The $\PSp^\pm(4,\bR)$ case was proved by Ryan Hoban in his thesis \cite{Hoban} using a similar invariant, with a small gap in the proof for the signature $(1,1)$ case.

\section{Preliminaries}

Let $V$ be a real vector space of even dimension $2n$. Denote by $\Gr(n,2n)$ the Grassmannian of $n$-dimensional subspaces in $V$. Let $U,W\in \Gr(n,2n)$ be transverse subspaces.  The projection to $U$ according to the splitting $V=U\oplus W$ will be denoted by $\pi_U^{W}$.

\begin{defn}
  The \emph{reflection} in the pair $U,W$ is the map
  \[\Refl{U}{W} = \pi_U^{W} - \pi_{W}^U.\]
\end{defn}

Given a splitting $V=U\oplus W$ into half-dimensional subspaces and a linear map $f:U \rightarrow W$, we will denote by $\graph(f)$ the subspace $\{u+f(u) ~|~ u\in U\}$. Any $n$-dimensional subspace which is transverse to $W$ is the graph of a unique linear map in this way. The subspace $\graph(f)$ is transverse to $U$ if and only if $f$ is invertible, and in this case $\graph(f)=\graph(f^{-1})$.

The main invariant that we use to prove rigidity is a generalized cross-ratio.

\begin{defn}
Let $U_1,U_2,U_3,U_4 \in \Gr(n,2n)$ such that $U_1,U_2$ and $U_3,U_4$ are transverse pairs. The \emph{cross-ratio} $[U_1,U_2;U_3,U_4]$ is the $\GL(U_1)$-conjugacy class of the endomorphism of $U_1$ defined by $\pi_{U_1}^{U_2}\circ \pi_{U_3}^{U_4}$.
\end{defn}

\begin{prop}
  Let $U_1,U_2,U_3,U_4 \in \Gr(n,2n)$ be pairwise transverse and write $U_2=\graph(f)$, $U_4=\graph(g)$ where $f,g\in \Hom(U_1,U_3)$. Then, the cross ratio $[U_1 , U_2 ; U_3 , U_4]$ is given by $f^{-1}\circ g$.
  \begin{proof}
    Let $v\in U_1$. Decomposing $v$ according to the splitting $V=U_3\oplus U_4$, there exists $u\in U_3$ and $u'\in U_1$ such that $v = u + (u' + g(u')) = u' + (u + g(u'))$.
    
    Therefore,  $u=-g(u')$ and $u'=v$, so $\pi_{U_3}^{U_4}(v)=u=-g(v)$.
    
    Similarly, for any $u\in U_3$ we have $\pi_{U_1}^{U_2}(u) = -f^{-1}(u)$. Hence,
    \[[U_1,U_2 ; U_3,U_4] = \pi_{U_1}^{U_2}\pi_{U_3}^{U_4}|_{U_1} = f^{-1}\circ g.\]
  \end{proof}
\end{prop}

\begin{nota}
  When comparing endomorphisms $f_1,f_2$ of different vector spaces $U_1,U_2$, we will use the notation $f_1\sim f_2$ to mean that there exists an isomorphism $g : U_1 \rightarrow U_2$ such that $f_2 = g\circ f_1\circ g^{-1}$.
\end{nota}

An elementary consequence of the previous proposition is that this generalized cross-ratio has the following symmetry whenever $U_1,U_2,U_3,U_4$ are pairwise transverse :
\begin{itemize}
    \item $[gU_1,gU_2;gU_3,gU_4] = g[U_1,U_2;U_3,U_4]g^{-1}$;
    \item $[U_1,U_2;U_4,U_3] = I-[U_1,U_2;U_3,U_4]$;
    \item $[U_1,U_4;U_3,U_2] = [U_1,U_2;U_3,U_4]^{-1}$;
    \item $[U_3,U_2;U_1,U_4] \sim [U_1,U_2;U_3,U_4]^{-1}$.
\end{itemize}

% \begin{prop}
% Let $L_1,L_2,L_3,L_4$ be pairwise transverse Lagrangians in $\bR^{2n}$, $n \geq 3$. Then, the invariants $\M(L_1,L_2,L_3)$ and $C=[L_1,L_2;L_3,L_4]$ characterize the quadruple up to the $\Sp(V,\omega)$ action.
% \begin{proof}
% We can normalize so that $L_1=\langle e_1,\dots, e_n\rangle$, $L_3=\langle e_{n+1},\dots,e_{2n}\rangle$, $L_2=\graph(dP)$, $L_4=\graph(dQ)$, where $P,Q$ are quadratic forms on $L_1$. The pairwise transversality assumption implies that $P,Q$ are nondegenerate and that there is no isotropic vector for both. Hence, by theorem \ref{thm:diagonalbasis} there is a basis in which both are diagonal.

% In this situation, the cross ratio $C$ is given by $PQ^{-1}$, and rescaling the basis so that $Q=I_{p,q}$ we get $P=I_{p,q}C$.
% \end{proof}
% \end{prop}

\section{Elements of finite order}
Since each reflection we consider is uniquely determined by a pair of transverse subspaces of dimension $n$, the cross ratio provides invariants of the linear transformations obtained by composing two reflections. Let us describe the cross-ratios of elements of finite order which can arise in this way.

\begin{prop}\label{prop:charpolIdentity}
The characteristic polynomial $p_R(\lambda)$ of the composition $T=\Refl{U_1}{U_2}\circ \Refl{U_3}{U_4}$ is related to the characteristic polynomial $p_C(\lambda)$ of the cross ratio
\[C=[U_1,U_2;U_3,U_4]\]
by the following equation :
\[p_T(\lambda) = (-4\lambda)^{n} p_C\left(\frac{(\lambda+1)^2}{4\lambda}\right).\]

%$\lambda^2 + 2(1-2\mu)\lambda + 1=0$. In fact, if $C$ has $n$ distinct eigenvalues, then $R$ has $2n$ distinct eigenvalues which are all the solutions of $\lambda^2 + 2(1-2\mu_i)\lambda +1=0$.
\begin{proof}
Consider a basis $e_1,\dots,e_{2n}$ of $V$ such that
\[U_1 = \langle e_1\dots e_n \rangle,\]
\[U_2 = \graph(f),\]
\[U_3 = \langle e_{n+1} \dots e_{2n}\rangle,\]
\[U_4 = \graph(g),\]
where $f,g : U_1\rightarrow U_3$ are linear maps. We will denote by $A,B$ the respective matrix expressions of $f,g$ with respect to the bases $e_1\dots,e_n$ and $e_{n+1}\dots e_{2n}$.

With this normalization, the cross ratio $C$ has matrix expression $A^{-1}B$, and the reflections satisfy
\[\Refl{U_1}{U_2} = \begin{pmatrix}
I & -2A^{-1}\\0 & -I
\end{pmatrix}\]
\[\Refl{U_3}{U_4} = \begin{pmatrix}
-I & 0\\
-2B & I \end{pmatrix}.\]
Their composition is therefore
\[T=\Refl{U_1}{U_2}\circ \Refl{U_3}{U_4} = \begin{pmatrix} 4A^{-1}B - I & -2A^{-1}\\
2B & -I \end{pmatrix}.\]
Using the block determinant formula
\[\det\begin{pmatrix}X & Y\\Z & W\end{pmatrix} = \det(XW-YZ),\]
which is valid whenever $Z$ and $W$ commute, we get
\begin{align*}
    \det(T-\lambda I) &= \det((\lambda+1)^2 I -4\lambda A^{-1}B)\\
                      &= (-4\lambda)^n\det\left(\frac{-(\lambda+1)^2}{4\lambda}I+A^{-1}B\right).
\end{align*}
% and so
% \[\mu = \frac{-(\lambda+1)^2}{4\lambda}\]
% is an eigenvalue of $C$, which we can rewrite as
% \[\lambda^2 + 2(2\mu-1)\lambda + 1 = 0.\]
% Similarly, if $\mu$ is an eigenvalue of $C$, from the equality of characteristic polynomials above, there must be an eigenvalue $\lambda$ of $R$ satisfying this quadratic equation.
\end{proof}
\end{prop}

\begin{prop} \label{prop:finiteorderEigenvalues}
Let $U_1,U_2,U_3,U_4$ be pairwise transverse $n$-dimensional subspaces. If the composition of reflections $T=\Refl{U_1}{U_2}\circ\Refl{U_3}{U_4}$ has order $N$ with $N\geq2$ in $\PGL(2n,\bR)$, then each eigenvalue $\mu_i$ of the cross-ratio $C=[U_1,U_2;U_3,U_4]$ is of the form $\mu_i = \sin^2(k_i\pi/(2N))$ for $k_i\in \bZ$.
\begin{proof}
If $T^n=\pm I$, all eigenvalues of $T$ must be $N$th roots of $1$ or $-1$. Then, by Proposition \ref{prop:charpolIdentity}, any eigenvalue $\mu$ of $C$ must satisfy
\[\lambda^2 + 2(2\mu-1)\lambda + 1 = 0\]
for some $\lambda$ an $N$-th root of $\pm 1$.

The solutions to the quadratic equation above are
\[\lambda = (2\mu-1) \pm 2\sqrt{\mu^2-\mu}.\]
The only possible solutions for $\lambda$ a root of unity satisfy $0\leq\mu\leq 1$. Changing variables to $\mu=\sin^2(\theta)$, we get $\lambda =-e^{\pm 2 i \theta}$ which is an $N$th root of $\pm 1$ only if $\theta=\frac{k\pi}{2N}$.
\end{proof}
\end{prop}

\section{Configurations of $6$-tuples}\label{sec:configurations}

Let $(U_1^+,U_1^-)$, $(U_2^+,U_2^-)$, $(U_3^+,U_3^-)$ be three pairs of half-dimensional subspaces in $V$ and denote by $R_i = \Refl{U_i^+}{U_i^-}$ the associated reflections. Assume that the subspaces $U_i^\pm$ are pairwise transverse. We can associate three cross-ratios to this configuration :
\[C_1 = [U_2^+,U_2^-;U_3^+,U_3^-],\]
\[C_2 = [U_1^+,U_1^-;U_3^+,U_3^-],\]
\[C_3 = [U_1^+,U_1^-;U_2^+,U_2^-].\]

\begin{defn}
The \emph{configuration space} of pairwise transverse $6$-tuples in $\Gr(n,2n)$ is
\[\conf^{(6)}(\Gr(n,2n)) := \left\{(U_1^+,U_1^-,U_2^+,U_2^-,U_3^+,U_3^-)\in \Gr(n,2n)^6 \right\}/\PGL(2n,\bR),\]
where $\PGL(2n,\bR)$ acts diagonally.

The three cross-ratios $C_1,C_2,C_3$ above define a map
\[\mathscr{C} : \conf^{(6)}(\Gr(n,2n)) \rightarrow \left(\GL(n,\bR)/\Inn(\GL(n,\bR))\right)^3\]
with values in triples of conjugacy classes of $\GL(n,\bR)$.
We will denote by
\[\conf^{(6)}_{C_1,C_2,C_3}(\Gr(n,2n)) := \mathscr{C}^{-1}(C_1,C_2,C_3)\]
the fiber of this map over the triple $(C_1,C_2,C_3)$.
\end{defn}

Assume that the reflections $R_1,R_2,R_3$ generate a $(k_1,k_2,k_3)$-triangle group. By proposition \ref{prop:finiteorderEigenvalues}, the set of conjugacy classes allowed for these cross-ratios in this case is finite and therefore those conjugacy classes must be fixed by any deformation of the triangle group.
Any continuous deformation of the group $\langle R_1,R_2,R_3 \rangle$ induces a path in the space $\conf^{(6)}_{C_1,C_2,C_3}(\Gr(n,2n))$.

The goal of this section will be to prove the following proposition :
\begin{prop}\label{prop:posDefFiniteM1}
  Let $\lambda_1,\lambda_2,\lambda_3 \in \bR$.
  If the polynomial
  \[\lambda_3(1-\lambda_2)x^2 + (\lambda_1+\lambda_2-\lambda_3-1)x + (1-\lambda_1)\]
  has distinct real roots, then the configuration space
  \[\conf^{(6)}_{\lambda_1 I,\lambda_2 I,\lambda_3 I}(\Gr(n,2n))\]
  is finite.
\end{prop}

Without loss of generality, we can assume the following:
\[U_1^+ = \langle e_{1}, \dots, e_{n}\rangle,\]
\[U_2^+ = \langle e_{n+1}, \dots, e_{2n}\rangle,\]
\[U_3^+ = \graph(N),\]
where $N: U_1^+ \rightarrow U_2^+$ is an invertible linear map. We can write the remaining three subspaces similarly as graphs $U_1^- = \graph(M_1)$, $U_2^-=\graph(M_2)$, and $U_3^-=\graph(M_3)$ of linear maps $M_i : U_1^+\rightarrow U_2^+$.

Then, The cross-ratios $C_i$ can be computed as functions of $M_1,M_2,M_3$ :
\begin{align}\label{eqn:crossratios}
C_1 &= (N-M_2)(N-M_3)^{-1}\\
C_2 &= M_1^{-1}(N-M_1)(N-M_3)^{-1}M_3 \nonumber\\
C_3 &= M_1^{-1}M_2 \nonumber
\end{align}

\begin{lem}\label{lem:polynomial_condition}
  If $C_1$,$C_2$,$C_3$ are scalars $C_i=\lambda_i I$, then $A_1 =N^{-1}M_1$ must satisfy the quadratic equation
  \[\lambda_3(1-\lambda_2)A_1^2 + (\lambda_1+\lambda_2-\lambda_3-1)A_1 + (1-\lambda_1)I = 0.\]
    \begin{proof}
    Assume $C_i = \lambda_i I$. Then, $M_2 = \lambda_3 M_1$. Substituting in the equation for $C_1$, we get
    \[\lambda_1(N - M_3) = N - \lambda_3 M_1\]
    which implies
    \[\lambda_1 M_3 = (\lambda_1-1) N  + \lambda_3 M_1.\]
    Finally, substituting in the equation for $C_2$ we find
    \[\lambda_2 I = M_1^{-1}(N-M_1)\left(N - \frac{\lambda_1-1}{\lambda_1}N - \frac{\lambda_3}{\lambda_1}M_1\right)^{-1}\left(\frac{\lambda_1-1}{\lambda_1}N + \frac{\lambda_3}{\lambda_1}M_1\right).\]
    This last equation simplifies to
    \[\lambda_3(1-\lambda_2)A_1^2 + (\lambda_1+\lambda_2-\lambda_3-1)A_1 + (1-\lambda_1)I = 0.\]
  \end{proof}
\end{lem}

Since the discriminant of the polynomial
 \[\lambda_3(1-\lambda_2)x^2 + (\lambda_1+\lambda_2-\lambda_3-1)x + (1-\lambda_1)\]
 
appearing in the previous proposition will be important in what follows, we will denote it by
\begin{align}
\phi(\lambda_1,\lambda_2,\lambda_3) = & 1 - 2(\lambda_1+\lambda_2+\lambda_3) + 2(\lambda_1\lambda_2 + \lambda_2\lambda_3 + \lambda_1\lambda_3) \\&+ \lambda_1^2+\lambda_2^2+\lambda_3^2 - 4 \lambda_1\lambda_2\lambda_3.\nonumber
\end{align}

%\begin{cor}\label{cor:posDefFiniteM1}
%   Let $U_1^\pm,U_2^\pm,U_3^\pm$ be six half-dimensional subspaces in $V$. If the cross-ratios $C_1,C_2,C_3$ are scalars $C_i=\lambda_i I$ and the polynomial
%   If the polynomial
%   \[\lambda_3(1-\lambda_2)x^2 + (\lambda_1+\lambda_2-\lambda_3-1)x + (1-\lambda_1)\]
%   has distinct real roots, then the configuration space
%   \[\conf^{(6)}_{\lambda_1 I,\lambda_2 I,\lambda_3 I}(\Gr(n,2n))\]
%   is finite.
  
  %$(U_1^+,U_1^-,U_2^+,U_2^-,U_3^+,U_3^-)$ is locally rigid.
  %\begin{proof}

   We now prove Proposition \ref{prop:posDefFiniteM1}. By Lemma \ref{lem:polynomial_condition}, with the normalizations of subspaces $U_i^\pm$ as above, the minimal polynomial of $A_1 = N^{-1}M_1$ is of degree at most two. By hypothesis, it has distinct real roots, and therefore $A_1$ is diagonalizable over $\bR$. This means that $N$ and $M_1$ are simultaneously diagonalizable by changing basis in both $U_1^+$ and $U_2^+$. The stabilizer of the pair $(U_1^+,U_2^+)$ in $\PGL(2n,\bR)$ acts on graphs precisely by simultaneous change of basis in $U_1^+$ and $U_2^+$. Therefore, we can assume that $N$ is given by the identity matrix and that $M_1$ is diagonal. The matrix equation
   \[\lambda_3(1-\lambda_2)M_1^2 + (\lambda_1+\lambda_2-\lambda_3-1)M_1 + (1-\lambda_1)I = 0\]
   then translates into $n$ real quadratic equations for the diagonal entries of $M_1$, each of which has $2$ solutions. The cross-ratio equations (\ref{eqn:crossratios}) then uniquely determine $M_2$ and $M_3$, finishing the proof of Proposition \ref{prop:posDefFiniteM1}.
   
%   For the second part, note that $I$ and $M_1$ are always simultaneously diagonalizable by the classical orthogonal diagonalization theorem and apply the same proof.
%  \end{proof}
%\end{cor}

\section{Triangle groups}
In this section, we first recall some facts about the geometric representation of a hyperbolic triangle group, and then prove the main local rigidity theorem for diagonally embedded representations into $\PGL(2n,\bR)$.

\subsection{Geometric representations in $\PGL(2,\bR)$}

The Lie group $\PGL(2,\bR)$ identifies naturally with the isometry group of the hyperbolic plane. The \emph{geometric representation} of a hyperbolic triangle group $\Delta(k_1,k_2,k_3)$ maps the generators to the three reflections in the sides of a triangle with interior angles $\frac{\pi}{k_i}$. It is uniquely defined up to conjugation, since hyperbolic triangles are determined up to isometry by their angles.

A different model for this geometric representation will be useful. Consider the \emph{Gram matrix}
\[M = \begin{pmatrix} 1 & -x & -y\\
-x & 1 & -z\\
-y & -z & 1\\ \end{pmatrix}.\]

Let $\Gamma$ be the reflection group generated by orthogonal reflections \[v\mapsto v -2 (v^tMe_i)e_i\] for the bilinear form given by $M$ in the canonical basis vectors $e_1,e_2,e_3$. The group $\Gamma$ preserves the bilinear form given by $M$. Assuming $x<1$, the determinant
\[\det(M) = 1 -x^2 -y^2 -z^2- 2xyz\]
is negative if and only if the signature of $M$ is $(2,1)$, in which case the group $\Gamma$ acts by hyperbolic isometries on the hyperboloid model of the hyperbolic plane.

The canonical or geometric representation of a triangle group is obtained when letting $x=\cos\left(\frac{\pi}{k_1}\right)$, $y=\cos\left(\frac{\pi}{k_2}\right)$, and $z=\cos\left(\frac{\pi}{k_3}\right)$. Hence, the polynomial $q(x,y,z) = 1 -x^2 -y^2 -z^2- 2xyz$ is negative whenever $x=\cos\left(\frac{\pi}{k_1}\right)$, $y=\cos\left(\frac{\pi}{k_2}\right)$, and $z=\cos\left(\frac{\pi}{k_3}\right)$ with $\frac{1}{k_1} + \frac{1}{k_2} + \frac{1}{k_3} < 1$.

The discriminant $\phi(\lambda_1,\lambda_2,\lambda_3)$ of the polynomial in Proposition \ref{prop:posDefFiniteM1} is related to the determinant of this Gram matrix by
\[\phi(\lambda_1,\lambda_2,\lambda_3) = - \frac{1}{4}q(1-2\lambda_1,1-2\lambda_2,1-2\lambda_3).\]

\begin{lem}
  The discriminant $\phi(\lambda_1,\lambda_2,\lambda_3)$ is strictly positive for $\lambda_i = \sin^2\left(\frac{\theta_i}{2k_i}\right)$ for $k_i$ satisfying $\frac{1}{k_1}+\frac{1}{k_2} + \frac{1}{k_3} < 1$.
  \begin{proof}
    By the observations above, we have
    \[q(1-2\lambda_1,1-2\lambda_2,1-2\lambda_3) = q\left(\cos\left(\frac{\pi}{k_1}\right),\cos\left(\frac{\pi}{k_2}\right),\cos\left(\frac{\pi}{k_3}\right)\right) < 0,\]
    and so
    \[\phi(\lambda_1,\lambda_2,\lambda_3) = -\frac{1}{4}q(1-2\lambda_1,1-2\lambda_2,1-2\lambda_3) > 0.\]
  \end{proof}
\end{lem}

For the geometric representation $\rgeom : \Delta(k_1,k_2,k_3)\rightarrow \PGL(2,\bR)$, we can use the normalization of the previous section with $n=1$. The matrices $M_i$ are scalars which we will denote by $m_i$ and the cross-ratios are real numbers which we denote by $c_i$.

\begin{figure}
    \centering
    \includegraphics[width=\textwidth]{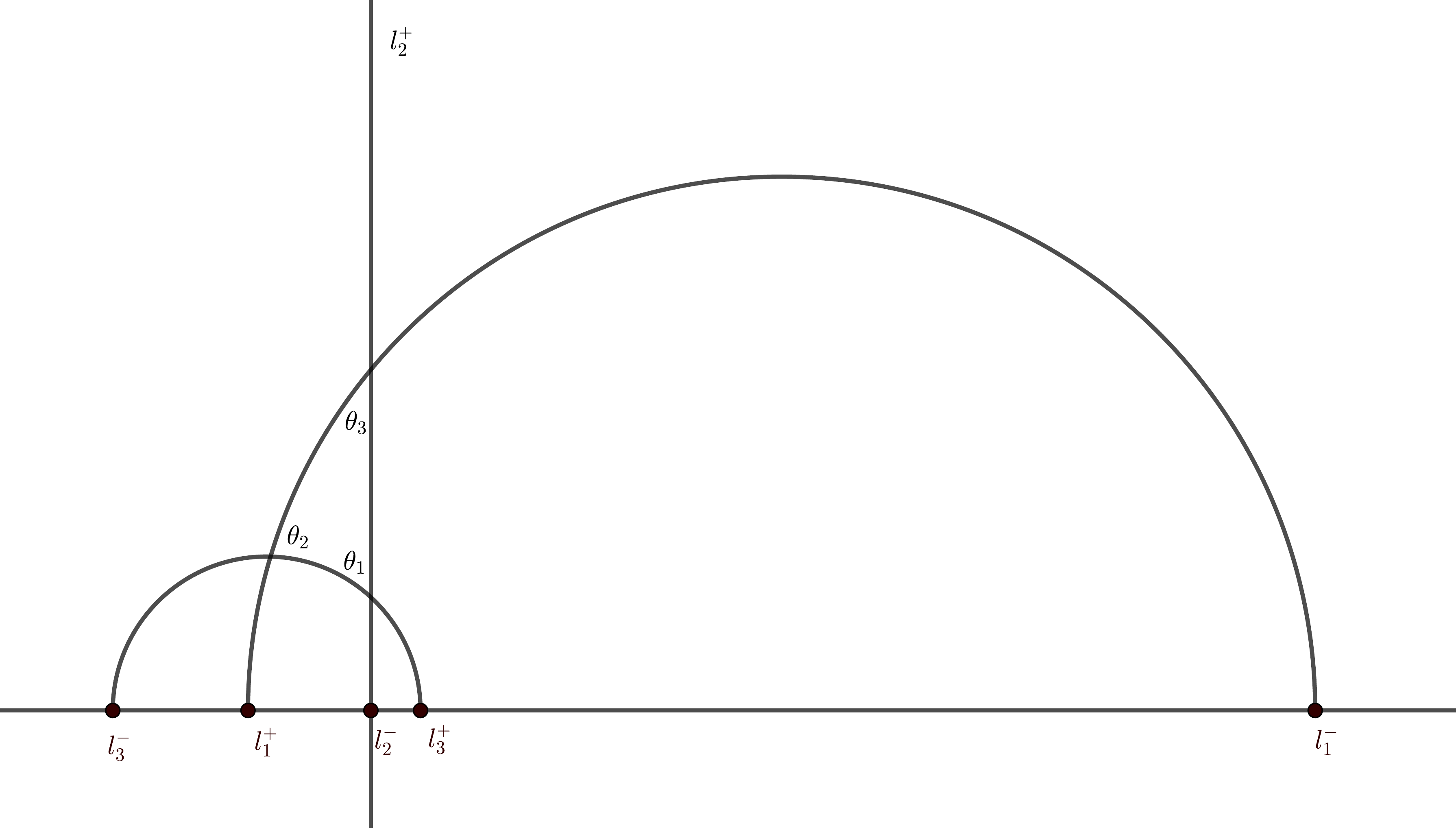}
    \caption{A hyperbolic triangle with labeled endpoints. The point $l_2^+$ is the point at infinity.}
    \label{fig:triangle}
\end{figure}

\begin{prop}\label{prop:crossratiosgeomsin}
  For the geometric representation of $\Delta(k_1,k_2,k_3)$, there is an ordering of the eigenspaces of the generators such that the cross-ratios have values $c_i=\sin^2(\frac{\pi}{2k_i})$.
  \begin{proof}
    Fix an orientation on the boundary of the hyperbolic plane and assume that the points $l_1^+,l_2^-,l_1^-,l_2^+ \in \partial \bH^2$ are placed in that order. Then, the angle between the two intersecting geodesics $\ell_1$, $\ell_2$ with respective endpoints $(l_1^+,l_1^-)$, $(l_2^+,l_2^-)$ in the hyperbolic plane is related to the cross-ratio of their endpoints by the formula
    \[\cos(\theta) = 1 - 2[l_1^+,l_1^- ; l_2^+,l_2^-].\]
    Note that the angle $\theta$ is always the angle given by the oriented arc between the boundary points $l_1^+$ and $l_2^-$ (or equivalently $l_1^-$ and $l_2^+$). Let $\Delta$ be a triangle in $\bH^2$ with interior angles $\theta_1 = \frac{\pi}{k_1}$, $\theta_2 = \frac{\pi}{k_2}$, $\theta_3 = \frac{\pi}{k_3}$. Order the endpoints $l_i^\pm$ of the geodesics extending the edges of $\Delta$ as in Figure \ref{fig:triangle}, so that
    \[[l_1^+,l_1^- ; l_2^+,l_2^-] = \frac{1-\cos(\theta_3)}{2}=\sin^2\left(\frac{\theta_3}{2}\right);\]
    \[[l_2^+,l_2^- ; l_3^+,l_3^-] = \frac{1-\cos(\theta_1)}{2}=\sin^2\left(\frac{\theta_1}{2}\right);\]
    \[[l_1^+,l_1^- ; l_3^+,l_3^-] = \frac{1-\cos(\theta_2)}{2}=\sin^2\left(\frac{\theta_2}{2}\right).\]
  \end{proof}
\end{prop}

\subsection{Diagonal representations into $\PGL(2n,\bR)$}\label{sec:rigidityPGL}
% \begin{thm}
%  If all crossratio invariants of a Lagrangian triangle group $\Delta$ are scalars, then $\Delta$ is a subgroup of a $(p,q)$-representation of $\Sp^\pm(2,\bR)$.
%  \begin{proof}
%   Denote by $L_1^\pm$, $L_2^\pm$, $L_3^\pm$ the three pairs of Lagrangians defining the reflections $R_1,R_2,R_3$, respectively. Since all crossratio invariants are scalars, by Proposition \ref{prop:scalarImpliesCircle} all of the $L_i^\pm$ lie on a common $(p,q)$-circle $S$ (WHY?), and so Corollary \ref{cor:commonCircleImpliesSubgroup} implies the theorem.
%  \end{proof}
% \end{thm}
Consider the diagonal embedding $i_{2n} : \PGL(2,\bR)\rightarrow \PGL(2n,\bR)$ which maps an element $A\in\PGL(2,\bR)$ to the diagonal action of $A$ on $(\bR^2)^n$. Let $\xi_{2n} : \RP^1 \rightarrow \Gr(n,2n)$ denote the map
\[\xi_{2n}(p) = \Span(\iota_1(p),\dots,\iota_n(p)),\]
where $\iota_j$ is the inclusion into the $j$th factor of $\bR^{2n} = \bR^{2} \oplus \dots \oplus \bR^2$. The map $\xi_{2n}$ is $i_{2n}$-equivariant.

\begin{lem}
  Let $\Refl{l_1}{l_2} \in \PGL(2,\bR)$ be the reflection in the pair of lines $l_1,l_2$. Then, $i_{2n}(\Refl{l_1}{l_2}) = \Refl{\xi(l_1)}{\xi(l_2)}$ is the reflection in their images by $\xi$.
  \begin{proof}
    Since the reflection $\Refl{l_1}{l_2}$ acts by the identity on $l_1$ and $-1$ on $l_2$, its image by $i_{2n}$ acts by the identity on $\xi_{2n}(l_1)$ and $-I$ on $\xi_{2n}(l_2)$.
  \end{proof}
\end{lem}

\begin{lem}\label{lem:crossratioscalar}
  Let $l_1,l_2,l_3,l_4 \in \RP^1$ be pairwise transverse with a given cross-ratio $[l_1,l_2;l_3,l_4] = c$. Then, the generalized cross-ratio satisfies
  \[[\xi_{2n}(l_1),\xi_{2n}(l_2);\xi_{2n}(l_3),\xi_{2n}(l_4)] = c I.\]
  \begin{proof}
    Fix a basis $e_1,e_2$ of $\bR^2$, and the corresponding basis $e_1^j = \iota_j(e_1)$, $e_2^j=\iota_j(e_2)$ for $j=1,\dots,n$. By transitivity of $\PGL(2,\bR)$ on pairs and equivariance of $\xi$, we can assume that $l_1=\Span(e_1)$ and $l_3=\Span(e_2)$. Then, we write $l_2=\graph(f)$ and $l_4=\graph(g)$ so that $c=\frac{g}{f}$.
    
    This implies that $\xi_{2n}(l_1) = \Span(e_1^1,\dots,e_1^n)$ and $\xi_{2n}(l_3) = \Span(e_2^1,\dots,e_2^n)$. Moreover, $\xi_{2n}(l_2) = \graph(F)$ and $\xi_{2n}(l_4) = \graph(G)$ where
    \[F(v_1 \oplus \dots \oplus v_n)= f(v_1) \oplus \dots \oplus f(v_n)\]
    and 
    \[G(v_1 \oplus \dots \oplus v_n)= g(v_1) \oplus \dots \oplus g(v_n).\]
    We find
    \[F^{-1}\circ G (v_1 \oplus \dots \oplus v_n) = \frac{g}{f} (v_1 \oplus \dots \oplus v_n) = c I\]
    and hence
    \[[\xi_{2n}(l_1),\xi_{2n}(l_2);\xi_{2n}(l_3),\xi_{2n}(l_4)] = F^{-1} G = c I.\]
  \end{proof}
\end{lem}

\begin{thm}\label{thm:rigidityPGL}
  Let $\Gamma = \Delta(k_1,k_2,k_3)$ be a hyperbolic triangle group. Then, the representations obtained by composing the geometric representation with the diagonal representation 
  \[i_{2n} : \PGL(2,\bR)\rightarrow \PGL(2n,\bR)\]
  are locally rigid.
  \begin{proof}
    Consider the configuration the six subspaces $U_1^\pm,U_2^\pm,U_3^\pm$ corresponding to eigenspaces of the generators. By Lemma \ref{lem:crossratioscalar} and Proposition \ref{prop:crossratiosgeomsin}, the three cross-ratios
    \[C_1 = [U_2^+,U_2^-;U_3^+,U_3^-],\]
    \[C_2 = [U_1^+,U_1^-;U_3^+,U_3^-],\]
    \[C_3 = [U_1^+,U_1^-;U_2^+,U_2^-],\]
    satisfy $C_i = \sin^2(\frac{\pi}{2k_i}) I$.
    
    Let $\rho_t : \Gamma \rightarrow \PGL(2n,\bR)$ be a continous family of representations such that $\rho_0 = i_{2n}\circ \rgeom$. Denote by $U_i^\pm(t)$ the corresponding paths of eigenspaces, and by $C_i(t)$ the cross-ratios as above.
    
    The finite order of each composition of generators $\rho(ab),\rho(bc),\rho(ac)$ is fixed, and since there are only finitely many conjugacy classes of elements of a fixed order in $\PGL(2n,\bR)$ these conjugacy classes must remain fixed throughout the deformation $\rho_t$. Therefore, by Proposition \ref{prop:finiteorderEigenvalues} the cross-ratios $C_1(t),C_2(t),C_3(t)$ must stay constant. By Proposition \ref{prop:posDefFiniteM1} the configuration $U_1^\pm(t),U_2^\pm(t),U_3^\pm(t)$ remains fixed up to the diagonal action of $\PGL(2n,\bR)$, and therefore $\rho_t$ is conjugate to $i_{2n}\circ \rgeom$.
  \end{proof}
\end{thm}

The \emph{expected dimension}, as defined in \cite{ALFHitchin}, is a heuristic dimension count for connected components of an orbifold surface character variety. Let $\Gamma=\Delta(k_1,k_2,k_3)$ and let $c_1,c_2,c_3$ be conjugacy classes of elements in $\PGL(n,\bR)$ which have order $k_1,k_2,k_3$ respectively. Then, the expected dimension of the component of the character variety which maps $(ab),(bc),(ac)$ to the respective conjugacy classes $c_1,c_2,c_3$ is defined to be
\[\dim_e(\chi_{c_1,c_2,c_3}(\Gamma,\PGL(n,\bR))) := \frac{1}{2}\left( \sum_{i=1}^3 \dim(c_i)  -2\dim \PGL(n,\bR) \right).\]
For a diagonal representations into $\PGL(2n,\bR)$, the conjugacy classes of the images of $(ab),(bc),(ac)$ have two eigenvalues $\zeta,\zeta^{-1}$, each with multiplicity $n$. Therefore, the centralizer in $\PGL(2n,\bR)$ of any element $M\in c_i$ is $(2n^2-1)$-dimensional, and we find
\begin{align*}
    \dim(c_i) &= \dim(\PGL(2n,\bR)) - \dim(Z_{\PGL(2n,\bR)}(M))\\
    &= (4n^2-1) - (2n^2-1)\\
    &= 2n^2.
\end{align*}
The expected dimension for the component containing diagonal representations is thus
\[\frac{1}{2}\left( 6n^2  -2\dim \PGL(n,\bR) \right) = \frac{1}{2}(6n^2 - 2(4n^2-1)) = -n^2+1,\]
which is negative unless $n=1$. Therefore, even though in this case the dimension count is not exact, it still accurately predicts local rigidity.

% \begin{thm}
%   Let $\Delta$ be a hyperbolic triangle group. Then, the Lagrangian reflection triangle groups obtained by composing the geometric representation with any signature $(p,q)$ diagonal representation 
%   \[i_{(p,q)} : \Sp^\pm(2,\bR)\rightarrow \Sp^\pm(2n,\bR)\]
%   are locally rigid.
% \end{thm}

\section{Diagonal representations into $\PSp^\pm(2n,\bR)$}
After recalling the relevant definitions for Lagrangian reflections in symplectic vector spaces, we will prove Theorem \ref{thm:sympl} using a method completely analogous to the $\PGL(2n,\bR)$ case.

Let $(V,\omega)$ be a real symplectic vector space of dimension $2n$ and let $\Sp(V,\omega)$ be its group of symplectic automorphisms. We will denote by $\Sp^\pm(V,\omega)$ the group of linear automorphisms $A$ of $V$ which satisfy $A^*\omega = \pm\omega$.

We first recall the theory of $\PSp^\pm(2n,\bR)$ orbits of transverse tuples of Lagrangians using the Maslov index. We roughly follow the exposition in Chapter 5 of \cite{ghyssignature}.

\begin{defn}
A \emph{Lagrangian subspace} $L\subset V$ is an $n$-dimensional subspace such that $\omega|_L=0$. The homogeneous space of all Lagrangians will be denoted by $\Lag(V)$ and called the Lagrangian Grassmannian.
\end{defn}

Let $L,L'$ be transverse Lagrangian subspaces of $V$. We will call the reflection $\Refl{L}{L'}$ a \emph{Lagrangian reflection}.

\begin{prop}
A Lagrangian reflection is anti-symplectic, that is, \[\left(\Refl{L}{L'}\right)^*\omega = -\omega.\]
In particular, $\Refl{L}{L'}\in \Sppm(V)$.
\begin{proof}
\begin{align*}
    \omega(\Refl{L}{L'}u,\Refl{L}{L'}v) &= \omega(\pi_L u - \pi_{L'} u, \pi_L v -\pi_{L'}v)\\
                              &= -\omega(\pi_L u,\pi_{L'}v) - \omega(\pi_{L'}u,\pi_L v)\\
                              &= -\omega(\pi_L u + \pi_{L'} u,\pi_L v + \pi_{L'} v)\\
                              &= -\omega(u,v).\qedhere
\end{align*}
\end{proof}
\end{prop}

Given any vector space $U$ of dimension $n$, the standard construction of the symplectic structure on the cotangent bundle of a manifold gives a symplectic vector space structure on $U\oplus U^*$. Explicitly, the symplectic form is given by
\[\omega_U(u \oplus \alpha,u' \oplus \alpha') := \alpha'(u) - \alpha(u').\]

Linear maps $f\in \Hom(U,U^*)$ parameterize $n$-dimensional subspaces of $U\oplus U^*$ which are transverse to $U^*$ via the graph construction $\graph(f) = \{ u \oplus f(u) ~|~ u\in U\}$. The subspace $\graph(f)$ is Lagrangian for $\omega_U$ if and only if the map $f$ is \emph{symmetric} in the following sense : $f(u)v = f(v)u$ for all $u,v\in U$. This is the standard notion of symmetric bilinear form on $U$ when identifying a bilinear form on $U$ with an element of $\Hom(U,U^*)$. The symmetric bilinear form is nondegenerate if and only if the corresponding graph is also transverse to $U$.

A pair of transverse Lagrangians $L_1,L_2$ in $V$ defines an isomorphism between the symplectic vector space $(V,\omega)$ and $(L_1\oplus L_1^*,\omega_{L_1})$. The isomorphism is given by $V\cong L_1\oplus L_2 \cong L_1\oplus L_1^*$ where the last isomorphism is the identity on $L_1$ and $v\rightarrow \omega(v,-)$ on $L_2$. To see that this is an isomorphism, it suffices to check
\[\omega_{L_1}(u\oplus \omega(v,\cdot), u'\oplus \omega(v',\cdot)) = \omega(v',u)-\omega(u',v) = \omega(u \oplus v,u' \oplus v').\]

\begin{defn}
Let $L_1,L_2,L_3$ be pairwise transverse Lagrangians in $V$. Identify $V\cong L_1\oplus L_3$ with $L_1\oplus L_1^*$ as above. The \emph{Maslov form} $B_{L_1,L_2,L_3}$ is the nondegenerate symmetric bilinear form on $L_1$ corresponding to $L_2$. The \emph{Maslov index} $\M(L_1,L_2,L_3)$ is the signature of $B_{L_1,L_2,L_3}$.
\end{defn}

The stabilizer in $\Sp(V,\omega)$ of the pair $L_1,L_3$ is isomorphic to $\GL(L_1)$ and its action on Lagrangians transverse to $L_3$ corresponds to the action on bilinear forms $B$ by the usual change of basis action $gB(u,v)= B(g^{-1}u,g^{-1}v)$.

The Maslov index is a complete $\Sp(V,\omega)$-invariant for triples of pairwise transverse Lagrangians. This means that any transverse triple of Lagrangians $L_1,L_2,L_3$ can be written as $L_1,L_2,\graph(I_{p,q})$, where $I_{p,q}$ is the standard diagonal matrix representing a bilinear form of signature $(p,q)$. The action of an anti-symplectic element in $g\in \Sp^\pm(V)$ reverses the signature, that is, if $\M(L_1,L_2,L_3) = (p,q)$, then $\M(gL_1,gL_2,gL_3) = (q,p)$.

\begin{defn}
Let $(p,q)$ be a a pair of natural numbers such with $p+q=n$. Let $L_1,L_2$ be a fixed pair of transverse Lagrangians in $(V,\omega)$ and $B$ a bilinear form of signature $(p,q)$ on $L_1$. The \emph{$(p,q)$-circle} defined by $B,L_1,L_2$ is the collection of Lagrangians which are multiples of $B$
\[\{\graph(\lambda B) \subset V ~|~ \lambda \in \bR\}\cup \{L_2\}.\]
\end{defn}

Identifying Lagrangians transverse to a fixed pair $L_1,L_2$ with nondegenerate bilinear forms, the problem of classifying quadruples of pairwise transverse Lagrangians translates to the problem of simultaneous diagonalization of bilinear forms. The following theorem tells us precisely when this is possible :

\begin{thm}[\cite{wonenburger}]\label{thm:diagonalbasis}
  Let $q_1,q_2$ be a pair of nondegenerate quadratic forms on a real vector space $V$. Denote by $\phi_1,\phi_2$ the isomorphisms $V\rightarrow V^*$ induced respectively by $q_1,q_2$. Then $q_1,q_2$ are simultaneously diagonalizable over $\bR$ if and only if the endomorphism $\phi_1^{-1}\phi_2$ is diagonalizable over $\bR$.
\end{thm}

% The quadruple invariant that we will use is a generalization of the cross-ratio of four points on the projective line. The same invariant was used in \cite{BurgerPozzetti} in order to study the geometry of $F$-tubes.

% \begin{defn}
% Let $L_1,L_2,L_3,L_4$ be Lagrangians such that $L_1,L_2$ and $L_3,L_4$ are transverse pairs. The \emph{cross-ratio} $[L_1,L_2;L_3,L_4]$ is the $\GL(L_1)$-conjugacy class of the endomorphism of $L_1$ defined by $\pi_{L_1}^{L_2}\circ \pi_{L_3}^{L_4}$.
% \end{defn}

% \begin{prop}
%   Let $L_1,L_2,L_3,L_4$ be pairwise transverse and write $L_2=\graph(f)$, $L_4=\graph(g)$ where $f,g\in \Hom(L_1,L_3)$. Then, the cross ratio $[L_1 , L_2 ; L_3 , L_4]$ is given by $f^{-1}\circ g$.
%   \begin{proof}
%     Let $v\in L_1$. Decomposing $v$ according to the splitting $V=L_3\oplus L_4$, there exists $u\in L_3$ and $u'\in L_1$ such that $v = u + (u' + g(u')) = u' + (u + g(u'))$.
    
%     Therefore,  $u=-g(u')$ and $u'=v$, so $\pi_{L_3}^{L_4}(v)=u=-g(v)$.
    
%     Similarly, for any $u\in L_3$ we have $\pi_{L_1}^{L_2}(u) = -f^{-1}(u)$. Hence, $[L_1,L_2,L_3,L_4] = \pi_{L_1}^{L_2}\pi_{L_3}^{L_4} = f^{-1}\circ g$.
%   \end{proof}
% \end{prop}

We will now show how to construct ``diagonal" homomorphisms $\PSL(2,\bR)\rightarrow \PSp(V)$ which preserve a $(p,q)$ circle for any signature $(p,q)$. These will be the analogs of the diagonal embedding $i_{2n}$ in the $\PGL(2n,\bR)$ case. Note that $\PSL(2,\bR)\cong \PSp(2,\bR)$ and $\PGL(2,\bR)\cong \PSp^\pm(2,\bR)$.

Let $(U,\omega)$ be a $2$-dimensional real symplectic vector space. Let $(W,b)$ be an $n$-dimensional real vector space equipped with a nondegenerate bilinear form $b$ of signature $(p,q)$. The vector space $U\otimes W$ is a symplectic vector space when equipped with the form $\omega \otimes b$.

\begin{prop}
  Denote
  \[L_u := \{u \otimes w ~|~ w \in W\}.\]
  Then, $L_u$ is a Lagrangian subspace only depending on the span of $u$ and the collection $S=\{L_u ~ |~ u\in U\}$ is a $(p,q)$-circle in $\Lag(U\otimes W)$.
\end{prop}

Now, since $\Aut(\omega)\cong \SL(2,\bR)$ and $\Aut(b)\cong O(p,q)$, the above construction defines a homomorphism
\[i_{(p,q)} : \SL(2,\bR) \times \Ort(p,q) \rightarrow \Sp(U\otimes W)\]
which preserves the $(p,q)$-circle $S$. More precisely, $\SL(2,\bR)$ acts on this circle equivariantly with respect to its action on $\mathbb{P}(U)$ and $\Ort(p,q)$ acts trivially.

\begin{prop}
  Three Lagrangians $L_1,L_2,L_3\subset V$ are contained in a unique $(p,q)$-circle, where $(p,q)=\M(L_1,L_2,L_3)$.
  \begin{proof}
    % Let $(W,b)=(L_1,b_{L_1,L_2,L_3})$ and $U = \{aL_1 + bL_3 | a,b\in\bR\}$, where $aL_1+bL_3$ is a formal linear combination of $L_1$ and $L_3$. Equip $U$ with the symplectic form 
    % \[\omega(aL_1 + bL_3,cL_1 + dL_3) = \frac 1 2 (a d - b c).\]
    % Then, $V\cong U\otimes W$ as symplectic vector spaces and the $(p,q)$-circle $S = \{L_u\}$ as above contains the three Lagrangians $L_i$.
    
    Let $S'$ be another $(p,q)$-circle containing $L_1,L_2,L_3$ and $L\in S'-S$. Then, writing $L_2=\graph(f)$, we must have $L=\graph(\lambda f)$ for some $\lambda$, which implies $L\in S$, a contradiction.
  \end{proof}
\end{prop}

\begin{prop}
  The homomorphism $i_{(p,q)}$ extends to a homomorphism of the double covers $\Sp^\pm(2,\bR) \rightarrow  \Sp^\pm(V)$. The orientation reversing elements in $\Sp^\pm(2,\bR)$ map to Lagrangian reflections.
  \begin{proof}
    The extension of the homomorphism is clear. For the second part, let $R\in \Sp^\pm(2,\bR)$ be a reflection. There exists a decomposition $U=\ell \oplus \ell'$ in which $R=1\oplus(-1)$. This implies that $i_{(p,q)}(R)|_{L_\ell} = 1$ and $i_{(p,q)}(R)|_{L_{\ell'}}=-1$. Since $L_\ell$ and $L_{\ell'}$ are complementary subspaces, this proves the claim.
  \end{proof}
\end{prop}

Finally, $i_{(p,q)}$ passes down to the quotients by the respective centers and induces a homomorphism, which we also denote $i_{(p,q)} : \PSp^\pm(2,\bR) \rightarrow \PSp^\pm(2n,\bR)$.

\begin{cor}\label{cor:commonCircleImpliesSubgroup}
  Let $L_1,L_2,L_3,L_4$ be Lagrangians lying on a common $(p,q)$-circle $S$. Then, $\Refl{L_1}{L_2}\Refl{L_3}{L_4}$ preserves $S$ and lies in the image of the corresponding $(p,q)$-representation of $\Sp^\pm(2,\bR)$.
\end{cor}

% \begin{prop}
%   Any two Lagrangians on a $(p,q)$ circle are transverse.
% \end{prop}

% \begin{prop}\label{prop:scalarImpliesCircle}
% The cross-ratio $[L_1,L_2;L_3,L_4]$ is a real multiple of the identity if and only if the $L_i$ lie on a common $(p,q)$-circle, and the signature $(p,q)$ is determined by the Maslov index $\M(L_1,L_2,L_3)$.
% \end{prop}

As in Section \ref{sec:configurations}, denote by $\conf^{(6)}_{C_1,C_2,C_3}(\Lag(2n))$ the space of pairwise transverse $6$-tuples of Lagrangians $L_1^\pm,L_2^\pm,L_3^\pm$ such that

\[C_1 = [L_2^+,L_2^-;L_3^+,L_3^-],\]
\[C_2 = [L_1^+,L_1^-;L_3^+,L_3^-],\]
\[C_3 = [L_1^+,L_1^-;L_2^+,L_2^-].\]

The analog of Proposition \ref{prop:posDefFiniteM1} for the symplectic case still holds :
\begin{prop}\label{prop:confRigidSymplectic}
  Let $\lambda_1,\lambda_2,\lambda_3\in \bR$. If the polynomial
  \[\lambda_3(1-\lambda_2)x^2 + (\lambda_1+\lambda_2-\lambda_3-1)x + (1-\lambda_1)\]
  has distinct real roots, then the configuration space
  \[\conf^{(6)}_{\lambda_1 I,\lambda_2 I,\lambda_3 I}(\Lag(2n))\]
  is finite.
  \begin{proof}
    The proof is a slight variation of the $\PGL(2n,\bR)$ case presented in Section \ref{sec:configurations}.
    Let $L_1^\pm, L_2^\pm, L_3^\pm$ be six pairwise transverse Lagrangian subspaces such that
    \[[L_2^+,L_2^-;L_3^+,L_3^-] = \lambda_1 I,\]
    \[[L_1^+,L_1^-;L_3^+,L_3^-] = \lambda_2 I,\]
    and
    \[[L_1^+,L_1^-;L_2^+,L_2^-] = \lambda_3 I.\]
    We take as our fixed pair of Lagrangians $L_1^+,L_2^+$ and write all the others as graphs of symmetric linear maps from $L_1^+$ to $L_2^+$, so
    \[L_3^+ = \graph(N)\]
    and
    \[L_i^- = \graph(M_i).\]
    By Lemma \ref{lem:polynomial_condition}, we know that the map $A_1 = N^{-1} M_1$ must satisfy the polynomial equation
    \begin{equation}\label{eqn:polsympl}
        \lambda_3(1-\lambda_2)A_1^2 + (\lambda_1+\lambda_2-\lambda_3-1)A_1 + (1-\lambda_1).
    \end{equation}
    In particular, the map $N^{-1}M_1$ is diagonalizable and therefore $N,M_1$ are simultaneously diagonalizable by Theorem \ref{thm:diagonalbasis}. We can therefore, by applying an element of $\stab(L_1^+,L_2^+)$, assume that $N,M_1$ are diagonal and so there are finitely many solutions to (\ref{eqn:polsympl}). Hence, there are only finitely many possibilities for $M_1,M_2,M_3$, finishing the proof.
  \end{proof}
\end{prop}

The orthogonal diagonalization theorem for quadratic forms allows us to prove a version of Proposition \ref{prop:confRigidSymplectic} with a hypothesis on the Maslov index instead of the roots of a polynomial.

\begin{prop}
  Let $(L_1^\pm,L_2^\pm,L_3^\pm)$ be a $6$-tuple of pairwise transverse Lagrangians such that $\M(L_1^+,L_2^+,L_3^+)=(n,0)$. Then, this tuple is isolated in the configuration space $\conf^{(6)}_{\lambda_1 I,\lambda_2 I,\lambda_3 I}(\Lag(2n))$.
  \begin{proof}
    The statement follows from the same argument as Proposition \ref{prop:confRigidSymplectic}, since if $\M(L_1^+,L_2^+,L_3^+)=(n,0)$, we can normalize $N$ to the identity matrix and orthogonally diagonalize $M_1$ by applying an element of $\stab(L_1^+,L_2^+,L_3^+)\cong \Ort(n)$.
  \end{proof}
\end{prop}

We deduce the local rigidity theorems for $\PSp^\pm(2n,\bR)$ (Theorem \ref{thm:sympl} and Theorem \ref{thm:symplpositive}) from the two previous propositions, using the same proof as in the $\PGL(2n,\bR)$ case (Theorem \ref{thm:rigidityPGL}).
% \begin{thm}
%   Let $\Gamma = \Delta(k_1,k_2,k_3)$ be a hyperbolic triangle group. Then, the representations obtained by composing the geometric representation with any signature $(p,q)$ diagonal representation
%   \[i_{(p,q)} : \PGL(2,\bR)\rightarrow \PSp^\pm(2n,\bR)\]
%   are locally rigid.
  
%   Similarly, 
% \end{thm}

Theorem \ref{thm:sympl} gives $\lceil \frac{n+1}{2} \rceil$ isolated points in the character variety $\chi(\Gamma,\PSp^\pm(2n,\bR))$ since this is the number of possible signatures $(p,q)$ with $p+q=n$, up to reversing $p$ and $q$.

% A similar count to that at the end of Section \ref{sec:rigidityPGL} gives that the expected dimension of the component of the character variety $\chi(\Gamma,\PSp^\pm(2n,\bR))$ containing diagonal $(p,q)$ representations is
% \[\dim_e(\chi_{c_1,c_2,c_3}(\Gamma,\PSp^\pm(2n,\bR)) = \]

\bibliographystyle{alpha}
\bibliography{biblio.bib}

\begin{thebibliography}{{Wie}18}

\bibitem[ALS18]{ALFHitchin}
Daniele {Alessandrini}, Gye-Seon {Lee}, and Florent {Schaffhauser}.
\newblock {Hitchin components for orbifolds}.
\newblock {\em arXiv e-prints}, page arXiv:1811.05366, Nov 2018.

\bibitem[BIW03]{BIW}
Marc Burger, Alessandra Iozzi, and Anna Wienhard.
\newblock Surface group representations with maximal {T}oledo invariant.
\newblock {\em C. R. Math. Acad. Sci. Paris}, 336(5):387--390, 2003.

\bibitem[GR16]{ghyssignature}
\'{E}tienne Ghys and Andrew Ranicki.
\newblock Signatures in algebra, topology and dynamics.
\newblock In {\em Six papers on signatures, braids and {S}eifert surfaces},
  volume~30 of {\em Ensaios Mat.}, pages 1--173. Soc. Brasil. Mat., Rio de
  Janeiro, 2016.

\bibitem[GW10]{gwcomponents}
Olivier Guichard and Anna Wienhard.
\newblock Topological invariants of {A}nosov representations.
\newblock {\em J. Topol.}, 3(3):578--642, 2010.

\bibitem[Hob09]{Hoban}
Ryan Hoban.
\newblock {\em Local rigidity of triangle groups in {S}p(4, {R})}.
\newblock ProQuest LLC, Ann Arbor, MI, 2009.
\newblock Thesis (Ph.D.)--University of Maryland, College Park.

\bibitem[LT18]{longth}
D.~Long and M.~Thistlethwaite.
\newblock The dimension of the hitchin component for triangle groups.
\newblock {\em Geometriae Dedicata (to appear)}, 2018.

\bibitem[Wei]{weir}
Elise~A. Weir.
\newblock The restricted hitchin component for triangle groups.
\newblock {\em in preparation}.

\bibitem[{Wie}18]{HTTSurvey}
Anna {Wienhard}.
\newblock {An invitation to higher Teichm\"uller theory}.
\newblock {\em arXiv e-prints}, page arXiv:1803.06870, Mar 2018.

\bibitem[Won66]{wonenburger}
Mar\'{\i}a~J. Wonenburger.
\newblock Simultaneous diagonalization of symmetric bilinear forms.
\newblock {\em J. Math. Mech.}, 15:617--622, 1966.

\end{thebibliography}

\end{document}